\documentclass[12pt]{amsart}
\usepackage{amsmath, amssymb, natbib}

\DeclareMathAlphabet{\curly}{U}{rsfs}{m}{n}

\newtheorem{thm}{Theorem}

\newtheorem{lem}{Lemma}[section]

\newcommand{\tS}{\widetilde{S}}
\newcommand{\tX}{\widetilde{X}}

\newcommand{\tR}{R}

\newcommand{\del}{\ensuremath{\delta}}
\newcommand{\lam}{\ensuremath{\lambda}}

\newcommand{\PPP}{\ensuremath{\mathbf{P}}}

\newcommand{\fl}[1]{{\ensuremath{\left\lfloor {#1} \right\rfloor}}}

\newcommand{\pfrac}[2]{{\left(\frac{#1}{#2}\right)}}

\newcommand{\be}{\begin{equation}}
\newcommand{\ee}{\end{equation}}
\newcommand{\benn}{\begin{equation*}}   
\newcommand{\eenn}{\end{equation*}}

\renewcommand{\(}{\left(}
\renewcommand{\)}{\right)}

\numberwithin{equation}{section}

\newif\ifdraft

\begin{document}

\title{Sharp probability estimates for generalized Smirnov statistics}

\author{Kevin Ford}

\dedicatory{Dedicated to the memory of Walter Philipp}

\address{Department of Mathematics, 1409 West Green Street, University
of Illinois at Urbana-Champaign, Urbana, IL 61801, USA}
\email{ford@math.uiuc.edu}

\date{\today}
\thanks{Research supported by National Science Foundation grants
DMS-0301083 and DMS-0555367.}

\begin{abstract} We give sharp, uniform estimates for the probability that
the empirical distribution function for $n$ uniform-$[0,1]$ random variables
stays to one side of a given line.
\end{abstract}

\maketitle

%
%
%
\section{Introduction}\label{sec:intro}
%
%
%

Let $U_1, \ldots, U_n$ be independent, uniformly distributed random
variables in $[0,1]$ and let $u>0$, $v>0$.
Our goal is to estimate
$$
Q_n(u,v) = \PPP\( F_n(t) \le \frac{vt+u}{n} \;\; (0\le t\le 1) \), 
$$
where $F_n(t)=\frac{1}{n}\sum_{U_i\le t} 1$ is the associated
 empirical distribution function. 
In 1939, N. V. Smirnov introduced the statistic
 $D_n^+ = \sqrt{n} \sup_{0\le t\le 1} (F_n(t)-t)$ and
proved \cite{Sm1} for each fixed $\lam\ge 0$ the asymptotic formula
\be\label{Smirnov}
\PPP(D_n^+\le \lam) = Q_n(\lam \sqrt{n},n) \to 1 - e^{-2\lam^2} 
\qquad (n\to\infty).
\ee
When $\lam_0 \le \lam = O(n^{1/6})$ with \emph{fixed} $\lam_0>0$,
sharper forms of \eqref{Smirnov} have been proven by a number of people
(e.g. \cite{Lau}; see also Ch. 9 of \cite{SW}), in particular
\be\label{Smirnov2}
\PPP(D_n^+\le \lam) =  1 - e^{-2\lam^2} \( 1 - \frac{2\lam}{3n^{1/2}} 
+ O\pfrac{\lam^4+1}{n} \).
\ee
Here and throughout the Landau $O-$symbol has its usual meaning: 
$f(\cdot)=O(g(\cdot))$
means $|f| \le c g$ for some constant $c$, which is independent of the
inputs to the function $f$.  Also, $f\ll g$ means $f=O(g)$ and $f
\asymp g$ means $f=O(g)$ and $g=O(f)$.

One may ask about the behavior of $Q_n(u,v)$ for a wider range of the
variables $u,v$.
The strong Koml\'os-Major-Tusn\'ady theorem \cite{KMT} implies
$$
|F_n(t)-t-n^{-1/2} B_n(t)| \ll \frac{\log n}{n} \qquad (0\le t\le 1)
$$
with probability $\ge 1 - O(1/n)$, where $B_n(t)$ is a Brownian bridge
process.  The order $\frac{\log n}{n}$ on the right side is also best possible
\cite{KMT} (see also Ch. 4 of \cite{CR}).  Since
$$
\PPP \( \sup_{0\le t\le 1} (B_n(t) - (at+b)) \le 0 \) = 1 - e^{-2b(a+b)},
$$
and writing 
$$
w=u+v-n,
$$
the KMT theorem implies the uniform estimate
\be\label{KMTQ}
\begin{split}
Q_n(u,v) &= O\pfrac{1}{n} + 1 - e^{-\frac{2(u+O(\log n))(w+O(\log n))}{n}}\\
&= 1 - e^{-2uw/n} + O\pfrac{(u+w+\log n)\log n}{n}.
\end{split}
\ee
This gives an asymptotic for $Q_n(u,v)$ provided $\frac{u}{\log n} \to\infty$,
$\frac{w}{\log n} \to \infty$, and $u+w=o(n/\log n)$ as $n\to\infty$.
In the author's recent paper \cite{F} on the distribution of divisors 
of integers, sharper information was needed for 
very small $u$ and $w$.
That paper includes a short proof of the crude bound
$Q_n(u,v) \ll \frac{(u+1)(w+1)^2}{n}$ uniformly in $n\ge 1$, $u\ge 0$, and 
$w\ge 0$.  

By using different methods, we prove here new uniform estimates,
which essentially remove the logarithm terms from the right side
of \eqref{KMTQ}.

%
%

\begin{thm}\label{Q}  
Uniformly in $u>0$, $w>0$ and $n\ge 1$, we have
$$
Q_n(u,v) = 1 - e^{-2 uw/n} + O\pfrac{u + w}{n}.
$$
\end{thm}

In particular, if $u\to\infty$, $w\to\infty$ and $u+w=o(n)$ as 
$n\to\infty$, then
$$
\frac{Q_n(u,v)}{1 - e^{-\frac{2uw}{n}}} \to 1.
$$

%
\section{A random walk with a barrier}
%

Exact formulas for $Q_n(u,v)$ are known, which we record below.

\begin{lem}\label{Qlem}
Assume $n\ge 1$ and $v>0$.  Then
\begin{enumerate}
\item If $n-v< u \le 1$, then $Q_n(u,v)=\frac{w}{v} (1 + u/v)
^{n-1}$;
\item  If $n-v<u<n $ and $u\ge 0$, then
$$
Q_n(u,v) = 1 - \frac{w}{v^n} \sum_{u<j\le n} \binom{n}{j} (v+u-j)^{n-j-1}
  (j-u)^j.
$$
\end{enumerate}
\end{lem}

Formula (i) is due to H. E. Daniels
\cite{Daniels} and (ii) is due to R. Pyke \cite{Pyke}.
The case $v=n$ in (ii) was earlier proved by Smirnov \cite{Sm1}.
Starting with (ii), one may use a more complicated version of the
complex analytic method of Lauwerier \cite{Lau}
to prove Theorem \ref{Q}.  This was carried out in an early
version of the author's paper \cite{F}, 
a sketch of which may be found in \cite{kol}
(the English paper \cite{kol2} includes a sketch of the argument below).
We present below an elementary, probabilistic proof
of Theorem \ref{Q}.
Rather than work with
(ii), we reinterpret $Q_n(u,v)$ in terms of a random walk.

\begin{lem}\label{walk}
Let $X_1, \ldots, X_{n+1}$ be independent random variables, each with
density function $e^{x-1}$ if $x \le 1$ and 0 if $x > 1$.  Put $S_0=0$ and
$S_j=X_1+\cdots+X_j$ for $j\ge 1$.  Then
$$
Q_n(u,v) = \PPP \left[ \max_{0\le j\le n} S_j < u \bigg| S_{n+1}=n+1-v \right].
$$
\end{lem}

\begin{proof}
Let $Y_1,\cdots,Y_{n+1}$ be independent random variables with
exponential distribution, and let $W_k=Y_1+\cdots+Y_{k}$ for $1\le k\le n+1$. 
Let $\xi_1,\ldots,\xi_n$ be the order statistics of $U_1,\ldots,U_n$,
so that $Q_n(u,v)$ is the probability that $\xi_j \ge \frac{j-u}{v}$ for
every $j$.
By a well-known theorem of R\'enyi \cite{Re}, 
the vectors $(\xi_1,\ldots,\xi_n)$ and 
$(W_1/W_{n+1},\ldots,W_n/W_{n+1})$ have identical distributions.
Similarly, given that $W_{n+1}=v$, the probability density function of the
vector $(W_1/v,\ldots, W_n/v)$ is identically $n!$ on the set
$\{ (x_1,\ldots,x_n) : 0\le x_1\le \cdots \le x_n\le 1\}$.  Therefore,
$$
Q_n(u,v)= \PPP \left[ \min_{1\le i\le n} (W_i-i) \ge -u \; \bigg| \; 
W_{n+1} = v \right]. 
$$
Putting $X_i=1-Y_i$ completes the proof.
\end{proof}

The sequence $0,S_1,S_2,\ldots$ can be thought of as a recurrent 
random walk on the real line, with $Q_n(u,v)$ being the probability
that the walk does not cross a barrier at the point $u$ given that
it ends at the point $n+1-v$ after $n+1$ steps.
A similar quantity may be defined for a random walk with
the $X_i$ having a different distribution.  In the paper \cite{walks},
an analog of Theorem \ref{Q} is proven for a general walk whose steps
$X_i$ have a
continuous or lattice distribution, but valid in a more limited range of the
variables.  More specifically, under appropriate conditions on $X_i$,
we prove that
$$
\PPP \left[ \max_{0\le j\le n-1} S_j < y \bigg| S_{n} = y-z \right]
 = 1 - e^{-2yz/n} + O\pfrac{y+z+1}{n}
$$
uniformly for $0\le y\le c\sqrt{n}$, $0\le z\le c\sqrt{n}$ ($c$ being any
fixed constant).

Kolmogorov used a relation similar to that in Lemma \ref{walk} 
in his seminal 1933 paper
\cite{Kol33} on the distribution of the statistic
$$
D_n = \sqrt{n} \sup_{0\le t\le 1} |F_n(t)-t|.
$$
Specifically, let $\tX_1,\tX_2,\ldots,\tX_n$ be independent random
variables with discrete distribution 
$$
\PPP [ \tX_j=r-1 ] = \frac{e^{-1}}{r!} \qquad (r=0,1,2,\ldots)
$$
and let $\tS_j=\tX_1+\cdots+\tX_j$ for $j\ge 1$.  
Like the variables $X_i$ in Lemma \ref{walk}, each $\tX_i$ has mean 0
and variance 1.  Kolmogorov proved that for integers $u\ge 1$,
\begin{align*}
\PPP ( \sup_{0\le t\le 1} |F_n(t)-t| \le u/n ) 
&= \frac{n! e^n}{n^n} \PPP \( \max_{0\le
  j\le n-1} |\tS_j| < u, \tS_n=0 \) \\
&= \PPP \(  \max_{0\le j\le n-1} |\tS_j| < u \, \bigg| \tS_n=0 \).
\end{align*}
Small modifications to the proof yield, for {\it integers} $u\ge 1$ and for
$n\ge 2$, that
$$
Q_n(u,n) =  \PPP \(  \max_{0\le j\le n-1} \tS_j < u \, \bigg| \tS_n=0 \).
$$
When $v\ne n$, however, it does not seem feasible to express $Q_n(u,v)$
in terms of the variables $\tS_j$.

Let $f_n$ be the density function for $S_n$ ($n=1,2,\ldots$).
The Central Limit Theorem for densities (e.g., Theorem 1 in \S 46 of
\cite{GK}) implies that 
for large $n$ and $|x|\ll \sqrt{n}$, $f_n(x)\approx (2\pi n)^{-1/2}
e^{-x^2/2n}$. 
However, there are asymmetries in the distribution for $|x|>\sqrt{n}$.
We have
\be\label{fmx}
f_n(x) = \begin{cases} \frac{(n-x)^{n-1}}{e^{n-x} (n-1)!} & \; x\le n
  \\ 0 & \; x>n \end{cases}
\ee
which is easily proved by induction on $n$.  

\begin{lem}\label{lemalpha}
Let $n\ge 2$.  Then
\begin{enumerate}
\item $f_n(x)$ is unimodular in $x$, with a maximum value $f_n(1)$,
 and $f_n(1) \sim \frac{1}{\sqrt{2\pi n}}$;
\item For $x\ge 0$, $f_n(1+x) \le f_n(1-x)$;
\item For each real $z\ge 0$, there is a unique number $b=b(n,z)$ satisfying
$0\le b\le z$ and $f_n(1-z)=f_n(1+z-b)$.
\end{enumerate}
\end{lem}

\begin{proof}
Item (i) follows from 
\be\label{fnder}
f_n'(x) = \frac{1-x}{n-x} f_n(x) \qquad (x<n)
\ee
and Stirling's formula.  For (ii), suppose $0\le x < n-1$.  Then
$$
\frac{f_n(1+x)}{f_n(1-x)} = e^{2x} \(1-\frac{x}{n-1}\)^{n-1}
\(1+\frac{x}{n-1}\)^{-(n-1)} = \exp\left\{ -2 \sum_{j=1}^\infty
\frac{x^{2j+1}}{(2j+1)(n-1)^{2j}} \right\} \le 1.
$$
Item (iii) follows immediately from (i) and (ii).
\end{proof}

Using properties of $b(n,z)$, we will prove a sharper form of Theorem \ref{Q}.

\begin{thm}\label{thmR}
Suppose $n\ge 1$, $1\le u\le \frac{n}{10}$, $1 \le w \le \frac{n}{10}$
and let $b=b(n+1,w)$.  Then
$$
Q_n(u,v) = 1 - \( 1 - \frac{u(2w-b)}{(n-w+b)(n+w-u)}\)^{n} + O\( \(
\frac{u+w}{n} + \frac{u w^2}{n^2} \) e^{-\frac{uw}{n+w-u}} \).
$$
\end{thm}

%
\section{A recurrence formula}
%

Our principal tool for estimating $Q_n(u,v)$ is a recurrence formula
based on the reflection principle for random walks : For $y\ge 0$ and
$y\ge x$, a recurrent random walk of $n$ steps 
that crosses the point $y$ and ends
at the point $x$ is about as likely as a random walk which ends at
$2y-x$ after $n$ steps.  For convenience, define
\newcommand{\DD}{\mathbf{D}}
$$
\tR_n(x,y) = f_n(x) \PPP \left[ \max_{0\le j\le n-1} S_j < y \bigg| S_n=x 
\right] = \DD\left[ \max_{0\le j\le n-1} S_j < y, S_n = x\right],
$$
where the last expression stands for the density function 
$\frac{d}{dx} \PPP[T_{n-1}< y,
S_n \le x]$.  From the reflection principle we expect that 
$\tR_n(x,y) \approx f_n(x) - f_n(2y-x)$.  The next lemma gives a precise
measure of the accuracy of the reflection principle for our specific random
walk. 

\begin{lem}\label{lemB}
For a positive integer $n\ge 2$, real $y > 0$, real $x$, and real
$a\ge 1$,
\be\label{lemBeq}
\tR_n(x,y) = f_n(x) - f_n(y+a) + \int_0^1 \sum_{k=1}^{n-1}
\tR_k(y+\xi,y) \( f_{n-k}(a-\xi)-f_{n-k}(x-y-\xi) \)\, d\xi.
\ee
\end{lem}

\begin{proof}
Define $T_j=\max(S_0,\ldots,S_j)$.
Start with
$$
\tR_n(x,y) = f_n(x) - f_n(y+a) + f_n(y+a) - \DD[T_{n-1} \ge y, S_n=x].
$$
If $S_n=y+a$, then there is a unique $k$, $1\le k\le n-1$, so that
$T_{k-1} < y$ and $S_k \ge y$.  Thus,
\benn\begin{split}
f_n(y+a) & = \sum_{k=1}^{n-1} \DD[T_{k-1} < y, S_k \ge y, S_n=y+a] \\
&= \sum_{k=1}^{n-1}  \int_0^1 \DD[T_{k-1} < y, S_k=y+\xi, S_n=y+a]\,
   d\xi \\
&=  \sum_{k=1}^{n-1}  \int_0^1 \tR_k(y+\xi,y) f_{n-k}(a-\xi)\, d\xi.
\end{split}\eenn
Similarly,
\benn\begin{split}
\DD[T_{n-1}\ge y, S_n=x] &=  \sum_{k=1}^{n-1} \DD[T_{k-1} < y, S_k \ge y,
  S_n=x] \\
&= \sum_{k=1}^{n-1} \int_0^1 \tR_k(y+\xi,y) f_{n-k}(x-y-\xi)\, d\xi.
\end{split}\eenn
\end{proof}

In Lemma \ref{lemB}, choosing $a \approx y-x-b(n,y-x)$ should make
$| f_{n-k}(a-\xi)-f_{n-k}(x-y-\xi) |$ small for small $k$.  Also, we expect
$\tR_k(y+\xi,y)$ to be small, especially for large $k$,
so the integral-sum on the right of
\eqref{lemBeq} will be treated as an error term.

The same argument provides an analogous formula when the steps in the
random walk have an arbitrary distribution (see \cite{walks}).

We next give a crude estimate for $\tR_n(x,y)$ when $x\ge y$ which will be used
on the right side of \eqref{lemBeq}.  

\begin{lem}\label{lemC}
If $k\ge 1$, $y\ge 0$, and $0\le \mu\le 1$, then 
$\tR_k(y+\mu,y) \ll \frac{y+1}{k} f_k(y)$.
\end{lem}

\begin{proof}
Without loss of generality, suppose $k\ge 10$ and $0\le y\le
\frac{k}{10}$.
By Lemma \ref{lemalpha} (ii), when $1\le j\le k-1$,
$f_j(4-\xi) \le f_j(\mu-\xi)$.  By Lemma \ref{lemB} (with $a=4$
and $x=y+\mu$) and \eqref{fnder},
$$
\tR_k(y+\mu,y) \le f_k(y+\mu) - f_k(y+4) = \int_{y+\mu}^{y+4}\frac{t-1}{k-t}
f_k(t)\, dt \ll \frac{(y+1)f_k(y)}{k}.
$$
\end{proof}

%
\section{Estimates for $f_n(x)$}
%

\begin{lem}\label{lemA}
We have
\begin{enumerate}
\item If $n\ge 20$ and $0\le z\le \frac{n}{10}$, then $b(n,z)\le \frac{z}3$ and
  $b(n,z) = \frac{2 z^2}{3(n-1)} +  O\pfrac{z^3}{n^2}$;
\item If $n\ge 1$ and $|x| \le \frac{n}{3}$, then $n^{-1/2} e^{-x^2/n}
  \ll f_n(x) \ll n^{-1/2} e^{-x^2/3n}$;
\item If $1\le h \le H \le 10x^2$, then $f_h(x) h^{-2} \ll
  f_H(x) H^{-2}$;
\item If $1\le k\le n$, then $f_k(x) \ll (n/k)^{1/2} f_n(x)$.
\end{enumerate}
\end{lem}

\begin{proof}
First, writing $b=b(n,z)$, we have
$$
\(1 - \frac{2z-b}{n-1+z} \)^{n-1} = e^{-2z+b}.
$$
Under the hypotheses of (i), let $t=\frac{2z-b}{n-1+z}$, so that
$0 \le t \le \frac15$ by Lemma \ref{lemalpha} (iii).
Then
$$
\frac{z}{n-1} = - \frac{\log(1-t)+t}{t} = \frac{t}{2} +
\frac{t^2}{3} + \cdots
$$
which implies
$$
t = 2\pfrac{z}{n-1} - \frac83 \pfrac{z}{n-1}^2 + O\(\pfrac{z}{n-1}^3\).
$$
The asymptotic for $b$ follows.  Since $\frac{t}{2} + \frac{t^2}{3} +
\cdots \le \frac35 t$,  $b\le \frac{z}{3}$ and this proves (i).

Item (ii) is trivial when $n<100$.
When $n\ge 100$, \eqref{fmx} and Stirling's
formula give
$$
f_n(x) \asymp n^{-1/2} e^{x-1} \(1-\frac{x-1}{n-1}\)^{n-1} = n^{-1/2}
\exp \left[ - \frac{(x-1)^2}{n-1} \sum_{m=2}^\infty \frac{1}{m}
  \pfrac{x-1}{n-1}^{m-2} \right].
$$
Since $|\frac{x-1}{n-1}| \le 0.35$, the sum on $m$ is between
$\frac12$ and $\frac23$, which proves (ii).

Since $f_n(x)n^{-2} \asymp n^{-5/2}$ for $n\ge
(\frac{x-1}{10})^2$, it suffices to prove (iii) when $H\le
(\frac{x-1}{10})^2$.  For $1\le h \le (\frac{x-1}{10})^2$ and $h>x$,
$$
f_h(x) h^{-2} \asymp g(h) := h^{-5/2} e^{x-1} \(1 - \frac{x-1}{h-1}\)^{h-1}.
$$
We have
$$
\frac{d}{dh} \log g(h) = \frac{-5}{2h} + \frac{x-1}{h-x} - \log\( 1 +
\frac{x-1}{h-x} \) > 0,
$$
and (iii) follows.
If $|x|\le \sqrt{n}$, Lemma \ref{lemalpha} (i) and part (ii)
above imply $f_k(x)\ll k^{-1/2}$ and $f_n(x) \gg n^{-1/2}$.
When $|x|>\sqrt{n}$, applying (iii) gives 
$f_k(x) \ll f_n(x) \ll (n/k)^{1/2} f_n(x)$, proving (iv).
\end{proof}

A useful corollary of Lemma \ref{lemC} and Lemma \ref{lemA} (iv) is
\be\label{eqD}
\tR_k(u+\xi,u) \ll \frac{n^{1/2} u f_{n+1}(u)}{k^{3/2}} \qquad (1\le k\le
n+1, u\ge 1, 0\le \xi\le 1).
\ee

%
%

\begin{lem}\label{lemE}
Suppose $n\ge 100$, $1\le w\le \frac{n}{10}$, $b=b(n+1,w)$ and $0\le
\xi\le 1$.
\begin{enumerate}
\item If $w^{3/2} \le h \le n$, then 
$$
|f_h(1+w-b-\xi)-f_h(1-w-\xi)| \ll \( \frac{w}{h} +
\frac{w^3}{h^2} \) f_{h} (1-w).
$$
\item If $2\sqrt{n} \le w \le \frac{n}{10}$ and $1\le k\le
n-3w$, then $f_{n+1-k}(1+w-b-\xi)$ and $f_{n+1-k}(1-w-\xi)$
  are each
$$
= f_{n+1}(1-w) \exp \left\{ \sum_{j=n-k}^{n} \( \frac{1}{2j} \(1 -
\frac{w^2}{j} \) + O\(\frac{w^3}{j^3}\) \) + O\pfrac{w}{n} \right\}.
$$
\end{enumerate}
\end{lem}

\begin{proof}
Assume $w^{3/2} \le h \le n$ and write
$$
\frac{f_h(1+w-b-\xi)}{f_h(1-w-\xi)} = e^{2w-b} \( 1 -
\frac{2w-b}{h-1+w+\xi} \)^{h-1} = e^E,
$$
where, by Lemma \ref{lemA} (i),
 \benn\begin{split}
E &= (2w-b) \( 1- \frac{h-1}{h-1+w+\xi} - \frac12
\frac{(h-1)(2w-b)}{(h-1+w+\xi)^2} + O\pfrac{w^2}{h^2} \) \\
&= \frac{2w-b}{h-1+w+\xi} \( w+\xi - \frac{2w-b}{2} \( 1 -
\frac{w+\xi}{h-1+w+\xi} \)\) + O\pfrac{w^3}{h^2} \\
&\ll \frac{w}{h} + \frac{w^3}{h^2}.
\end{split}\eenn
By hypothesis, $E\ll 1$ and hence
 \benn\begin{split}
|f_h(1+w-b-\xi)-f_h(1-w-\xi)| &= f_h(1-w-\xi) |e^E-1| \ll |E|
f_h(1-w-\xi) \\
&\le |E| f_h(1-w) \ll \(  \frac{w}{h} + \frac{w^3}{h^2} \)
f_h(1-w). 
\end{split}\eenn
This proves (i).

To prove (ii), we write
\be\label{Ea}
\begin{split}
\frac{f_{n+1-k}(1+w-b-\xi)}{f_{n+1}(1+w-b)} &= \frac{f_{n+1}(1+w-b-\xi)}
{f_{n+1}(1+w-b)}
\prod_{j=n-k}^{n} \frac{f_{j}(1+w-b-\xi)}{f_{j+1}(1+w-b-\xi)} \\
&= e^{A+B_{n-k}+\cdots +B_{n}},
\end{split}
\ee
say.  By \eqref{fmx} and the hypothesis on $w$, $A \ll \frac{w}{n}$ and
\benn\begin{split}
B_j &= 1 + (j-1)\log\(1-\frac{1}{j-w+b+\xi}\) + \log\(1 +
  \frac{w-b-\xi}{j-w+b+\xi}\) \\
& = \frac{1}{j-w+b+\xi} \( 1 - \frac{j-1}{2(j-w+b+\xi)}
  - \frac{(w-b-\xi)^2}{2(j-w+b+\xi)} + O\( \frac{1}{j} + \frac{w^3}{j^2} 
  \)\) \\
& =  \frac{1}{j} \( \frac12 - \frac{w^2}{2j} \) + O\pfrac{w^3}{j^3}.
\end{split}\eenn
Arguing similarly,
\be\label{Eb}
\frac{f_{n+1-k}(1-w-\xi)}{f_{n+1}(1-w)}= e^{C+D_{n-k}+\cdots+D_{n}},
\ee
where $C\ll \frac{w}{n}$ and 
$D_j = \frac{1}{2j}(1-\frac{w^2}{j}) +  O(w^3/j^3).$
Combining \eqref{Ea}, \eqref{Eb}, the above estimates for $A$, $B_j$,
$C$ and $D_j$, and the relation $f_{n+1}(1-w)=f_{n+1}(1+w-b)$ 
concludes the proof of (ii).
\end{proof}

%
\section{proof Theorem \ref{thmR}}
%

Without loss of generality, suppose $n\ge n_0$, where $n_0$ is a large
absolute constant. 
We apply Lemma \ref{lemB} with $a=1+w-b$, where $b=b(n+1,w)$, obtaining
\be\label{thmR1}
\tR_{n+1}(n+1-v,u) = f_{n+1}(u+1-w) - f_{n+1}(u+1+w-b) + 
\sum_{k=1}^{n} \Delta_k,
\ee
where
$$
|\Delta_k| \le \max_{0\le \xi\le 1} \tR_k(u+\xi,u) \bigl|
f_{n+1-k}(1+w-b-\xi)-f_{n+1-k}(1-w-\xi) \bigr|.
$$
If $n\le u^2$, then $f_k(u)/k \ll f_{n+1}(u)/n$ for $1\le k\le n$ 
by Lemma \ref{lemA} (iii).  If $n>u^2$, then we have
$$
\frac{f_k(u)}{k} \ll \begin{cases} u^{-1} f_{\fl{u^2}}(u) \asymp 
u^{-3} & \text{ if }k\le u^2 \\ k^{-3/2} & \text{ if }k>u^2 \end{cases}
$$
by Lemma \ref{lemA} (ii), (iii).  In both cases,
\be\label{sumfky}
\sum_{k=1}^n \frac{f_k(u)}{k} \ll \( 1 + \frac{n^{1/2}}{u}\) f_{n+1}(u).
\ee

Suppose that $1\le w\le 2\sqrt{n}$, so that $b=O(1)$.
We will prove that
\be\label{thmRa}
\sum_{k=1}^{n} |\Delta_k| \ll \frac{u+w}{n} f_{n+1}(u) \ll
\frac{u+w}{n^{1/2}} f_{n+1}(u) f_{n+1}(1-w) \quad (1\le w\le 2\sqrt{n}).
\ee
The second inequality follows from the first and Lemma \ref{lemA}
(ii).  Let $h=n+1-k$, $h_0=\fl{w^{3/2}}$ and $h_1 = \fl{w^2/10}$.
Choose $n_0\ge 2^{10}$ so that $h_0\le n/2$.  For $1\le h\le h_0$,
\eqref{eqD} and Lemma \ref{lemA} (ii) give
$$
\Delta_{n+1-h} \ll \frac{uf_{n+1}(u)}{n} \max_{0\le \xi\le 1} \( 
f_h(1+w-b-\xi)+f_h(1-w-\xi) \) \ll \frac{uf_{n+1}(u)}{nw^3}.
$$
If $h_0<h\le h_1$, then \eqref{eqD}, Lemma \ref{lemE} (i) and Lemma
\ref{lemA} (ii),(iv) imply
$$
\Delta_{n+1-h} \ll \frac{uf_{n+1}(u)}{n} \( \frac{w}{h} +
\frac{w^3}{h^2} \) f_h(1-w) \ll \frac{uf_n(u)}{n} \( \frac{w}{h_1} +
\frac{w^3}{h_1^2} \) f_{h_1}(1-w) \ll \frac{uf_{n+1}(u)}{nw^2}. 
$$
When $h_1 < h \le \frac{n}{2}$, \eqref{eqD} and Lemma \ref{lemE} (i)
imply
$$
\Delta_{n+1-h} \ll \frac{uf_{n+1}(u)}{n} \frac{w}{h} f_h(1-w) \ll 
\frac{uw f_{n+1}(u)}{nw^{3/2}}.
$$
Summing on $h\le \frac{n}{2}$ we obtain
\be\label{thmRd}
\sum_{1\le h\le \frac{n}{2}} |\Delta_{n+1-h}| \ll \frac{uf_{n+1}(u)}{n}.
\ee
Lemma \ref{lemC}, Lemma \ref{lemE} (i) and \eqref{sumfky} imply
$$
\sum_{n/2 < h\le n} |\Delta_{n+1-h}| \ll \frac{uw}{n^{3/2}} \sum_{1\le k
  < n/2+1} \frac{f_k(u)}{k} \ll \frac{u+w}{n} f_{n+1}(u). 
$$
Combined with \eqref{thmRd}, this proves \eqref{thmRa}.

Next, suppose $2\sqrt{n} < w \le \frac{n}{10}$ and set
$$
K= \fl{ \min\( n-C_0w, \frac{n^3}{w^3} \)},
$$
where $C_0$ is a large absolute constant.
When $1\le k\le K$, apply Lemma \ref{lemC} and Lemma \ref{lemE} (ii),
observing that for each $j\le n$, $\frac{1}{2j} (1-w^2/j) \le 
-\frac{w^2}{3j^2} \le -\frac{w^2}{3n^2}$.  If $k\le n/2$, then
\be\label{Delk}
\Delta_k \ll u \frac{f_k(u)}{k} \( \frac{w}{n} +
\frac{kw^3}{n^3} \) e^{-kw^2/(10n^2)} f_{n+1}(1-w).
\ee
When $n/2 < k \le K$, 
$$
\Delta_k \ll \frac{u f_k(u) f_{n+1}(1-w)}{k} e^{-\frac{w^2}{6(n-k)}} \( 
\exp \left\{ C_1 \(\frac{w^3}{(n-k)^2} + \frac{w}{n} \) \right\} - 1 \) 
$$
for an absolute constant $C_1$.  If in addition $n-k \ge w^{3/2}$, then
$$
 e^{-\frac{w^2}{6(n-k)}} \( 
\exp \left\{ C_1 \(\frac{w^3}{(n-k)^2} + \frac{w}{n} \) \right\} - 1 \) 
\ll \( \frac{w^3}{(n-k)^2}+\frac{w}{n} \) e^{-\frac{w^2}{6(n-k)}},
$$
which implies \eqref{Delk}.  If $C_0 w \le n-k < w^{3/2}$ and we take
$C_0=20C_1$, then
$$
 e^{-\frac{w^2}{6(n-k)}} \( 
\exp \left\{ C_1 \(\frac{w^3}{(n-k)^2} + \frac{w}{n} \) \right\} - 1 \) 
\ll e^{-\frac{w^2}{12(n-k)}} \ll n^{-3} e^{-\frac{kw^2}{10n^2}},
$$
and \eqref{Delk} follows in this case as well.

By Lemma \ref{lemA} (iv), \eqref{sumfky} and \eqref{Delk},
\be\label{thmRh}\begin{split}
\sum_{k\le K} |\Delta_k| &\ll u f_{n+1}(1-w) \left[ \frac{w}{n}\sum_{k\le K}
  \frac{f_k(u)}{k} + \frac{w^3 f_{n+1}(u)}{n^{5/2}} \sum_{k=1}^\infty k^{-1/2}
  e^{-kw^2/(10n^2)} \right] \\
&\ll f_{n+1}(1-w) f_{n+1}(u) \( \frac{uw^2}{n^{3/2}} + \frac{w}{n^{1/2}} \).
\end{split}\ee
When $k>K$, we combine
Lemma \ref{lemA} (i), (iii) and Lemma \ref{lemE} (ii) to obtain
\begin{align*}
f_{n+1-k}(1+w-b-\xi)+f_{n+1-k}(1-w-\xi) &\ll f_{n+1-K}(1+w-b-\xi)+
f_{n+1-K}(1-w-\xi) \\
&\ll e^{-Kw^2/(10n^2)} f_{n+1}(1-w).
\end{align*}
Together with \eqref{eqD}, this gives
$$
\sum_{K<k\le n} |\Delta_k| \ll un^{1/2} f_{n+1}(u) f_{n+1}(1-w)
e^{-Kw^2/(10n^2)} \sum_{K<k\le n} \frac{1}{k^{3/2}}.
$$
If $2\sqrt{n} < w \le n^{2/3}$, then $K=\fl{n-3w}$ and 
$$
e^{-Kw^2/(10n^2)} \sum_{K<k\le n} \frac{1}{k^{3/2}}
\ll \frac{w}{n^{3/2}} e^{-w^2/(20n)} \ll \frac{w^2}{n^2}.
$$
If $n^{2/3} < w \le \frac{n}{10}$, then $K\ge n^3/2w^3$ and
$$
e^{-Kw^2/(10n^2)} \sum_{K<k\le n} \frac{1}{k^{3/2}} \ll
(n^3/w^3)^{-1/2} e^{-n/20w} \ll \frac{w^2}{n^2}.
$$
Therefore,
$$
\sum_{K<k\le n} |\Delta_k| \ll f_{n+1}(1-w) f_{n+1}(u) \frac{uw^2}{n^{3/2}}.
$$
Combined with \eqref{thmRh}, we have
\be\label{thmRi}
\sum_{k=1}^{n} |\Delta_k| \ll f_{n+1}(1-w) f_{n+1}(u) \( \frac{uw^2}{n^{3/2}}+
\frac{w}{n^{1/2}} \) \qquad (2\sqrt{n} < w \le n/10).
\ee

Combining \eqref{thmR1}, \eqref{thmRa} and \eqref{thmRi} with
Lemma \ref{walk}, in all cases we have
$$
Q_n(u,v) = 1 -\frac{f_{n+1}(u+1+w-b)}{f_{n+1}(u+1-w)} + O\(
\frac{n^{1/2}f_{n+1}(1-w)f_{n+1}(u)}{f_{n+1}(u+1-w)} \left[ \frac{u+w}{n} +
  \frac{uw^2}{n^2} \right] \).
$$
By the definition of $b$,
\begin{align*}
\frac{f_{n+1}(u+1+w-b)}{f_{n+1}(u+1-w)} &= \frac{f_{n+1}(u+1+w-b) 
f_{n+1}(1-w)}{f_{n+1}(1+w-b)
  f_{n+1}(u+1-w)} \\
&= \( 1 - \frac{u(2w-b)}{(n-w+b)(n+w-u)}\)^{n}.
\end{align*}
Also, by Stirling's formula,
\begin{align*}
\frac{n^{1/2}f_{n+1}(1-w)f_{n+1}(u)}{f_{n+1}(u+1-w)} &=
\frac{n^{1/2}(n+1)^n}{e^{n+1} n!}
\pfrac{(n+1-u)(n+w)}{(n+1)(n+w-u)}^{n} \\
&\ll \( 1 - \frac{u(w-1)}{(n+1)(n+w-u)}\)^{n} \ll e^{-\frac{uw}{n+w-u}},
\end{align*}
which concludes the proof of Theorem \ref{thmR}.
\qed

%
\section{proof Theorem \ref{Q}}
%

We may assume $0 \le u\le \del n$ and $0\le w\le
\del n$ for a small, fixed, positive $\del$.  
If $0\le u\le 1$ and $0\le w\le \del n$, Lemma
\ref{Qlem} (i) implies $Q_n(u,v) \ll w/n$.   When $0\le w\le 1$ and
$1\le u\le \del n$, Lemma \ref{lemC} implies $Q_n(u,v) \ll \frac{u}{n}.$
When $1\le u\le \del n$ and $1\le w\le \del n$, 
we may assume that $n$ is large. 
The error term in Theorem \ref{thmR} is
$$
\ll \frac{u+w}{n} + \frac{w}{n} \cdot \frac{uw}{n} e^{-\frac{uw}{2n}}
\ll \frac{u+w}{n}.
$$
When $uw>n^{4/3}$, the main terms are
$$
1 - O(e^{-\frac12 n^{1/3}}) = 1 - e^{-\frac{2uw}{n}} + O\pfrac{1}{n}.
$$
When $uw\le n^{4/3}$, the main terms are, by Lemma \ref{lemA} (i),
\begin{align*}
& = 1-\exp\left[ - \frac{u(2w-b)n}{(n-w+b)(n+w-u)} +
  O\pfrac{(uw)^2}{n^3}  \right] \\
&= 1- \exp\left[ -\frac{2uw}{n} \(1 +  O\pfrac{u+w}{n}\)\right]
  +O\pfrac{(uw)^2}{n^3}  \\
&= 1 - e^{-\frac{2uw}{n}} + O\( \frac{u+w+(uw)^{1/2}}{n} \) \\
&= 1 -  e^{-\frac{2uw}{n}} + O\pfrac{u+w}{n}.
\end{align*}
\qed

\medskip
{\bf Acknowledgments.}  The author expresses thanks to
Valery Nevzorov, Walter Philipp,
Steven Portnoy, and Jon Wellner for helpful conversations.
The author also thanks the referee for suggestions on improving the
exposition.

%
%

\bibliographystyle{plain}
\bibliography{smirnov}

\end{document}